\documentclass{article}
\usepackage[british]{babel}
\usepackage{amsmath,amssymb,graphicx,enumerate,xcolor}

\newcommand{\assign}{:=}
\newcommand{\tmcolor}[2]{{\color{#1}{#2}}}
\newcommand{\tmmathbf}[1]{\ensuremath{\boldsymbol{#1}}}
\newcommand{\tmop}[1]{\ensuremath{\operatorname{#1}}}
\newcommand{\tmtextbf}[1]{\text{{\bfseries{#1}}}}
\newcommand{\tmtextit}[1]{\text{{\itshape{#1}}}}
\newenvironment{enumerateroman}{\begin{enumerate}[i.] }{\end{enumerate}}
\newenvironment{itemizedot}{\begin{itemize} }{\end{itemize}}
\begin{document}

\title{An interactive visualisation for all $2 \times 2$ real matrices, with
applications to conveying the dynamics of iterative eigenvalue algorithms}

\author{Ran Gutin (rg120@ic.ac.uk, Imperial College London)}

\maketitle

\begin{abstract}
  We present two interactive visualisations of $2 \times 2$ real matrices,
  which we call v1 and v2. v1 is only valid for PSD matrices, and uses the
  spectral theorem in a trivial way -- we use it as a warm-up. By contrast, v2
  is valid for \tmtextit{all} $2 \times 2$ real matrices, and is based on the
  lesser known theory of Lie Sphere Geometry. We show that the dynamics of
  iterative eigenvalue algorithms can be illustrated using both. v2 has the
  advantage that it simultaneously depicts many properties of a matrix, all of
  which are relevant to the study of eigenvalue algorithms. Examples of the
  properties of a matrix that v2 can depict are its Jordan Normal Form and
  orthogonal similarity class, as well as whether it is triangular, symmetric
  or orthogonal. Despite its richness, using v2 interactively seems rather
  intuitive.
\end{abstract}

\section{Introduction}

We introduce two visualisations in this paper for use in self-study or
teaching. We call them v1 and v2. They are both very simple geometrically.
They are as follows:
\begin{itemizedot}
  \item v1 is perhaps obvious to those well-versed in linear algebra. We treat
  it is a warm-up. It is valid only for PSD matrices. A matrix is portrayed as
  an ellipse centred at the origin of $\mathbb{R }^2$. This method is an easy
  consequence of the spectral theorem.
  
  \item v2 is less obvious to those who know linear algebra because it uses
  machinery from Lie Sphere Geometry. It is valid \tmtextit{for all} $2 \times
  2$ real matrices. In ``half of cases'', a matrix is portrayed as an
  unordered pair, consisting of an \tmtextit{oriented circle}, and its
  \tmtextit{orientation-reversed inversion} through the unit circle.
\end{itemizedot}
v2 has the advantage that many properties of an arbitrary $2 \times 2$ matrix
$M$ can quickly be inferred from it:
\begin{itemizedot}
  \item Its Jordan Normal Form.
  
  \item Which matrices it's orthogonally similar to.
  
  \item Its condition number.
  
  \item Whether it's upper triangular, lower triangular, diagonal, symmetric,
  or orthogonal.
  
  \item A geometric construction for matrix-vector multiplication: $M v$.
  
  \item A geometric construction for matrix-matrix multiplication: $M K$.
\end{itemizedot}
While there are some ways to extend v1 to cover all $2 \times 2$ real matrices
(for instance by depicting the polar decomposition of $M$), we don't find that
most of the above properties hold. A visualisation that depicts the polar
decomposition might still be useful for depicting the dynamics of SVD-finding
or polar decomposition-finding algorithms. We don't consider those here as
they're rather obvious.

Our approach aims to convey the dynamics of the \tmtextit{unshifted} QR
algorithm, without using formulas or theory. Prior pedagogical approaches to
the QR algorithm have been considered before in the literature
{\cite{watkins2008qr,buurema1958geometric}}, but these are more formula-driven
and less visual than ours. Extensions of v2 to $3 \times 3$ matrices are
formally possible, but raise some issues which are not investigated here.

\section{Visualisation v1 (as a warm-up)}

The visualisation is \tmtextbf{valid only for PSD symmetric matrices}. A $2
\times 2$ PSD matrix $M$ is visualised by $\{ M v : v \in \mathbb{R}^2, \| v
\|_2 = 1 \}$. The visualisation of $M$ is thus always an ellipse, thanks to
the spectral theorem.

The \tmcolor{blue}{blue ellipse is the input}, the {\color[HTML]{006400}green
ellipse is one iteration of the LR algorithm}, and the \tmcolor{red}{red
ellipse is one iteration of the QR algorithm}.

\raisebox{0.0\height}{\includegraphics[width=3.13659320477502cm,height=2.92748917748918cm]{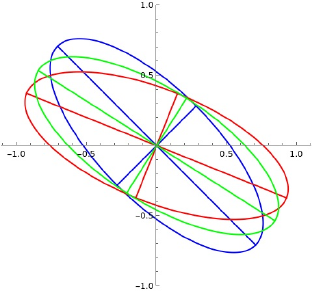}}\raisebox{0.0\height}{\includegraphics[width=3.13659320477502cm,height=2.92748917748918cm]{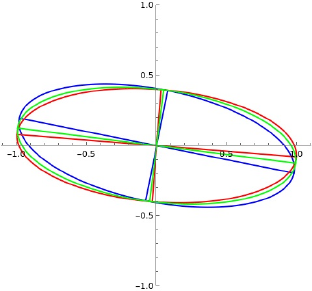}}\raisebox{0.0\height}{\includegraphics[width=3.13659320477502cm,height=2.92748917748918cm]{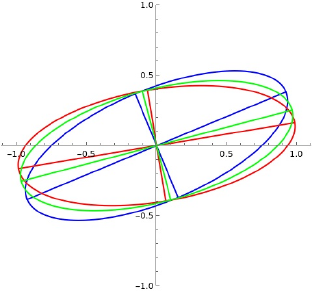}}\raisebox{0.0\height}{\includegraphics[width=3.13659320477502cm,height=2.92748917748918cm]{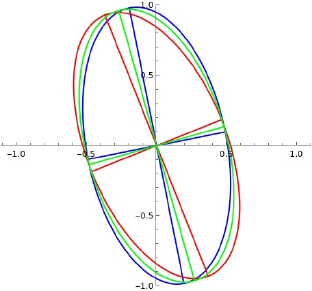}}

\raisebox{0.0\height}{\includegraphics[width=3.13659320477502cm,height=2.92748917748918cm]{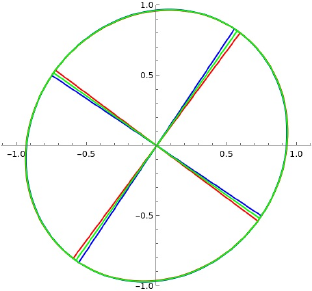}}\raisebox{0.0\height}{\includegraphics[width=3.13659320477502cm,height=2.92748917748918cm]{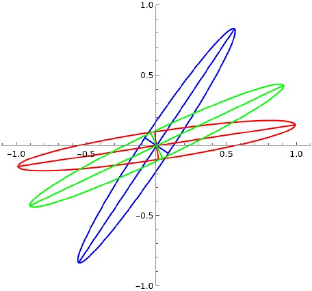}}

The important points are:
\begin{itemizedot}
  \item Convergence happens in every case, though not necessarily quickly.
  
  \item Ellipses where the major semi-axis is parallel to the $x$-axis are
  attractive fixed points. Convergence near here is fast.
  
  \item Ellipses where the major semi-axis is parallel to the $y$-axis are
  repulsive fixed points. Convergence near here is slow.
  
  \item When the eigenvalues are close to each other, convergence is slow. The
  ellipses look nearly like circles.
  
  \item When one eigenvalue is close to zero, and the other is not, then
  convergence is fast. The ellipses then look squeezed. They look ``opposite
  to circles''. \tmtextit{Shifting} (at least in the PSD case) can be seen as
  squeezing an ellipse.
\end{itemizedot}

\section{Visualisation v2}

\subsection{Oriented cycles, informally}

An example of an \tmtextit{oriented cycle} (note: not necessarily a
\tmtextit{circle}) and its \tmtextit{orientation-reversed inversion} is given
below. The pair of oriented cycles is in red. There is also a faint dark-blue
circle, which represents the unit circle. (Note that we have superimposed two
unit circles: One for each possible orientation).

\raisebox{0.0\height}{\includegraphics[width=5.9119277187459cm,height=6.21802439984258cm]{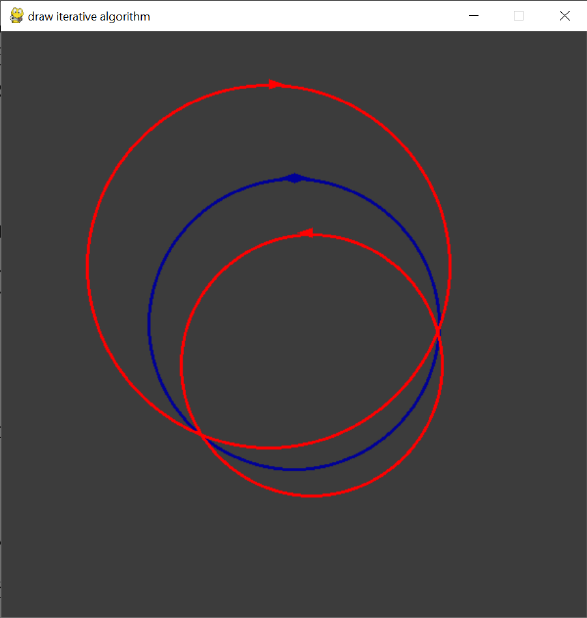}}

To explain the terminology: A \tmtextit{cycle} is either a point, a circle, or
a line. An \tmtextit{oriented cycle} is either:
\begin{itemizedot}
  \item a point $p$, which may be $\infty$,
  
  \item an \tmtextit{oriented non-point circle}, which is a pair $(B,
  B^{\ast})$, where $B$ is a disk, and $B^{\ast}$ is the set of points either
  inside or outside of $B$,
  
  \item an \tmtextit{oriented line}, which is a half-plane $H$. The underlying
  line is given by the boundary of $H$.
\end{itemizedot}
For each non-point cycle, there are at most two orientations. The orientation
can be depicted by drawing an arrow head.

An \tmtextit{inversion} through the unit circle (which is depicted in dark
blue) is effectively the map $z \mapsto 1 / \overline{z}$ over the complex
plane. If one thinks of an oriented cycle $C$ as a point moving in the complex
plane (with the orientation being the direction of motion, ignoring both the
exact position and speed of the moving point) then the \tmtextit{inversion of
C} is given by applying the inversion to the moving point.

An \tmtextit{orientation-reversed inversion} is a transformation which
performs an inversion, and then reverses all orientations. The order in which
these two operations are performed does not matter.

We denote the orientation-reversed inversion of an oriented cycle $C$ by
$\tmop{ORI} (C)$,

\subsection{Projective real line as a circle}

The \tmtextit{projective real line} $\mathbb{\tmop{RP}}^1$ is the set
$\mathbb{R} \cup \{ \infty \}$. The set $\mathbb{\tmop{RP}}^1$ can be depicted
as a circle. In our case, it is the unit circle in dark blue.

\raisebox{0.0\height}{\includegraphics[width=2.5974025974026cm,height=2.65861537452447cm]{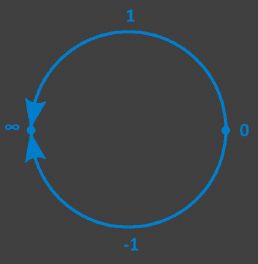}}

Instead of using the $\infty$ symbol, it is better to coordinatise its
elements using homogeneous coordinates: The homogeneous coordinate $[x : y]$
(where $x \in \mathbb{R}, y \in \mathbb{R},$and at least one of $x$ or $y$ is
non-zero) represents the element $\frac{x}{y}$ of $\mathbb{\tmop{RP}}^1$.
There is a canonical projection map $\pi : \mathbb{R}^2 \rightarrow
\mathbb{\tmop{RP}}^1, (x, y)^T \mapsto [x : y]$.

Every $2 \times 2$ real matrix $M$ canonically acts on $\mathbb{\tmop{RP}}^1$
by $M. [x : y] = \pi (M (x, y)^T)$. Two matrices $M$ and $K$ admit the same
canonical action on $\mathbb{\tmop{RP}}^1$ iff $M = \lambda K$ for $\lambda$
some non-zero scalar.

A \tmtextit{hyperbolic line} is a cycle orthogonal to the unit circle, i.e. it
meets the unit circle at a right angle. A \tmtextit{hyperbolic spear} is an
\tmtextit{oriented} hyperbolic line. In particular, a hyperbolic spear is
always its own orientation-reversed inversion. A given hyperbolic spear can be
denoted by the ordered pair $(p, q)$ where both $p$ and $q$ are in
$\mathbb{\tmop{RP}}^1$, and are the endpoints of the hyperbolic spear. Note
that the term \tmtextit{hyperbolic} here might appear confusing given that
there is no hyperbola, but this instead a reference to hyperbolic geometry,
which we will not dwell upon.\footnote{These are indeed the lines of
hyperbolic geometry, as depicted in the Poincare disk model of hyperbolic
geometry. The term ``spear'' usually refers to oriented lines in geometries in
which the term ``line'' can be given a sensible meaning.}

Given a non-singular $2 \times 2$ real matrix $M$, we can consider the set of
hyperbolic spears which represents its graph $\Gamma (M) \assign \{ (p, M.p) :
p \in \mathbb{\tmop{RP}}^1 \}$. We then attempt to somehow summarise $\Gamma
(M)$ using a single geometric figure.

\subsection{The non-negative determinant case}

We intend to show that any non-zero matrix with non-negative determinant can
be depicted as an oriented cycle and its orientation-reversed inversion. This
includes all non-zero singular matrices.

Given an invertible matrix $M$ with non-negative determinant, there is a
unique pair of oriented cycles $\{ C, \tmop{ORI} (C) \}$ such that each
hyperbolic spear in $\Gamma (M)$ is orientedly tangent to both $C$ and
$\tmop{ORI} (C)$. We will verify this claim in the appendix using Lie Sphere
Geometry.

We have an example of this below, where the elements of $\Gamma (M)$ are in
dark red, and both $C$ and $\tmop{ORI} (C)$ are in bright red:

\raisebox{0.0\height}{\includegraphics[width=5.9119277187459cm,height=6.21802439984258cm]{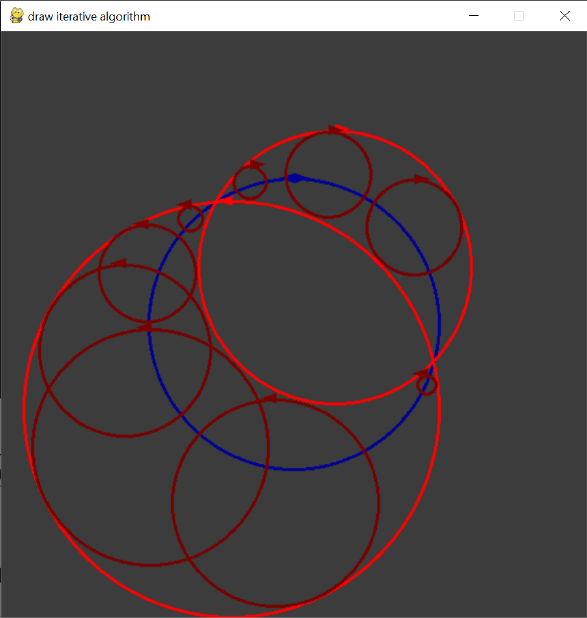}}\raisebox{0.0\height}{\includegraphics[width=5.9119277187459cm,height=6.21802439984258cm]{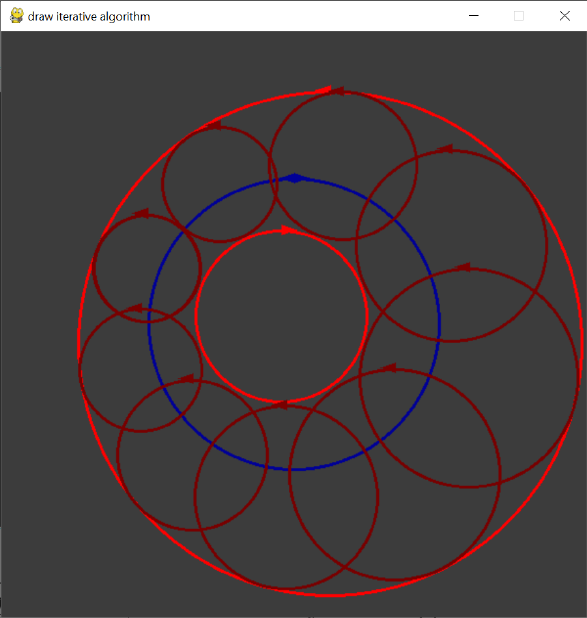}}

We may also pass to the limit where $M$ is singular (but $M \neq
\tmmathbf{0}$). The claim that there is a unique $\{ C, \tmop{ORI} (C) \}$
corresponding to $M$ still holds, as we will verify in the appendix.

\subsection{The non-positive determinant case}

Consider a negative determinant matrix $M$. Essentially, its representation is
a \tmtextit{continuously-oriented }\tmtextit{hyperbolic line} where the
orientation is in the closed interval $[- 1, + 1]$ instead of in the set $\{ -
1, + 1 \}$.

To understand why: Consider the special case where $M$ is diagonal.

\raisebox{0.0\height}{\includegraphics[width=5.9119277187459cm,height=6.21802439984258cm]{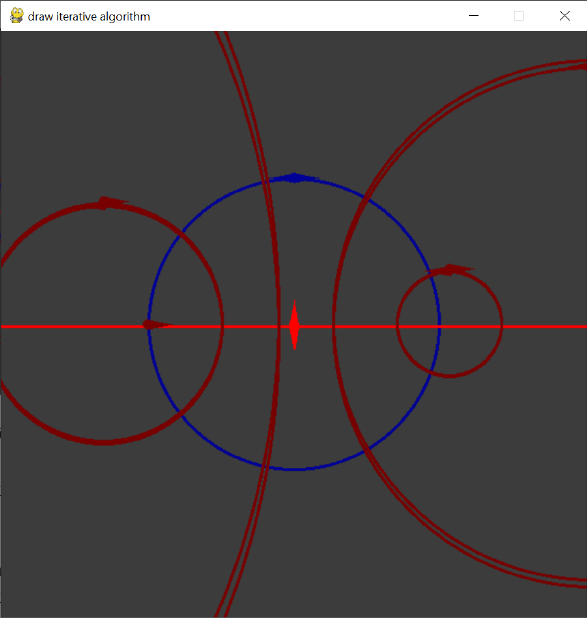}}\raisebox{0.0\height}{\includegraphics[width=5.9119277187459cm,height=6.21802439984258cm]{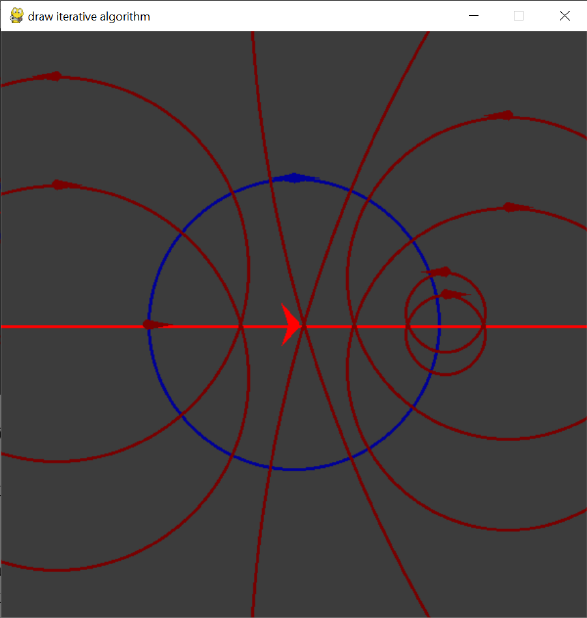}}\raisebox{0.0\height}{\includegraphics[width=5.9119277187459cm,height=6.21802439984258cm]{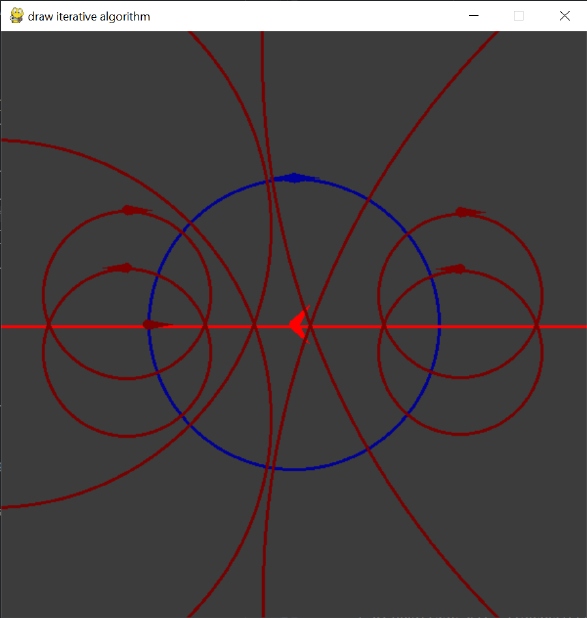}}

In dark red, we see $\Gamma (M)$. We see that all hyperbolic spears in $\Gamma
(M)$ meet the hyperbolic line in bright red at two set angles. The two set
angles are depicted by the \tmtextit{two} bright red arrow heads. We say that
the orientation of the bright red line, which represents our matrix, is
$\theta \in [- 1, 1]$, if the angles of approach are $\pm \frac{\pi}{2} (1 +
\theta)$. These are thus two different values, depicted using two arrow heads.
In the event that $\theta \in \{ - 1, 1 \}$, the two angles of approach will
become the same, and so the bright red arrow heads will point in the same
direction as each other. In that same event, the matrix $M$ will be singular,
and the visualisation will look identical to the non-negative determinant
case.

The case where $M$ is not diagonal looks as follows:

\raisebox{0.0\height}{\includegraphics[width=5.9119277187459cm,height=6.21802439984258cm]{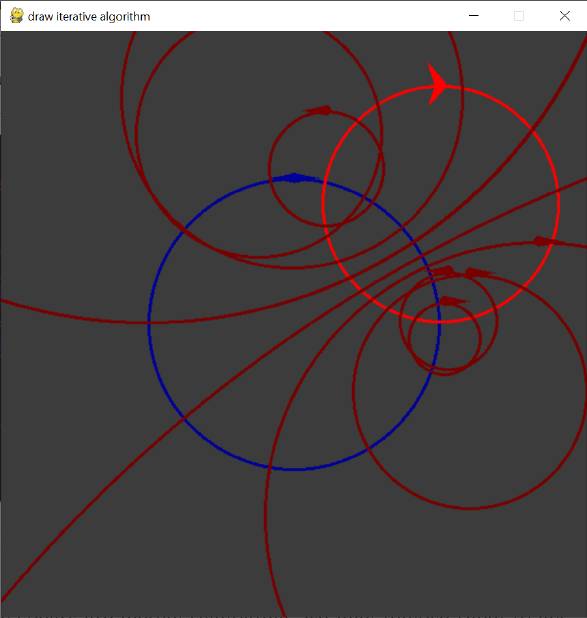}}

The only difference is that the bright red circle (which is the hyperbolic
line depicting $M$) is no longer a straight line.

Since the dark red hyperbolic spears are merely scaffolding, we may remove
them to get:

\raisebox{0.0\height}{\includegraphics[width=5.9119277187459cm,height=6.21802439984258cm]{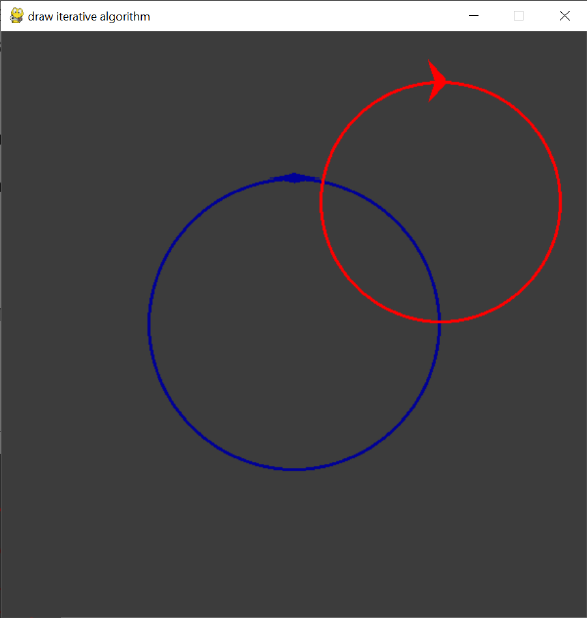}}

\subsection{Visualising the naive QR algorithm for all non-negative
determinant $2 \times 2$ real matrices}

The QR algorithm is complicated enough to benefit from this visualisation.
Here are some visualisations of the QR algorithm. The input matrix is given in
red. The first iteration is in light blue. Subsequent iterations get lighter
and lighter, until they're white. The user is expected to control one of the
red circles by making it:
\begin{itemizedot}
  \item Translate left, right, up, or down
  
  \item Shrink or expand
  
  \item Reverse in orientation
\end{itemizedot}
This determines the movement of the other red circle as well, which may then
become a straight line. In this way, the user can traverse all possible input
$2 \times 2$ real matrices with non-negative determinants.

\begin{table}[h]
  \begin{tabular}{lll}
    1 & 2 & 3\\
    \raisebox{0.0\height}{\includegraphics[width=5.9119277187459cm,height=6.21802439984258cm]{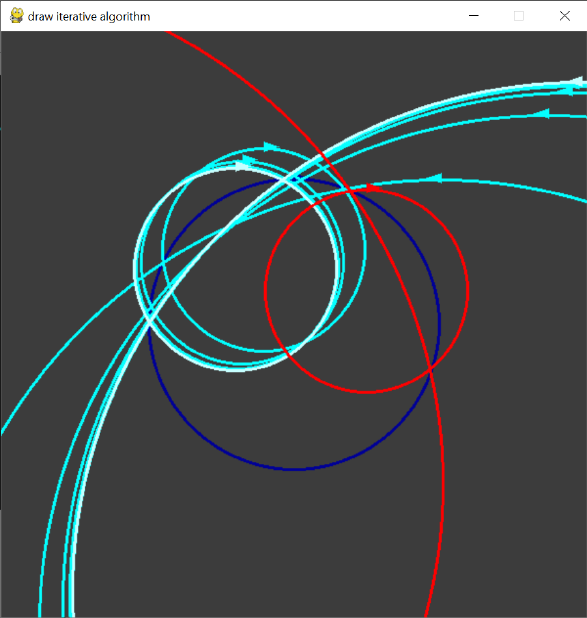}}
    &
    \raisebox{0.0\height}{\includegraphics[width=5.9119277187459cm,height=6.21802439984258cm]{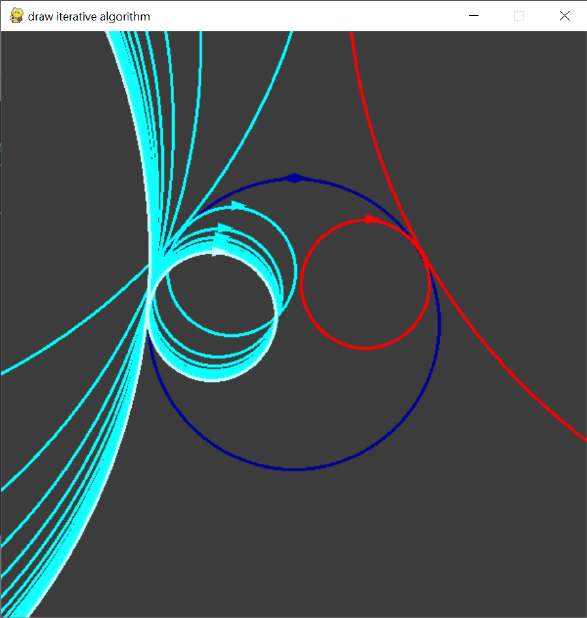}}
    &
    \raisebox{0.0\height}{\includegraphics[width=5.9119277187459cm,height=6.21802439984258cm]{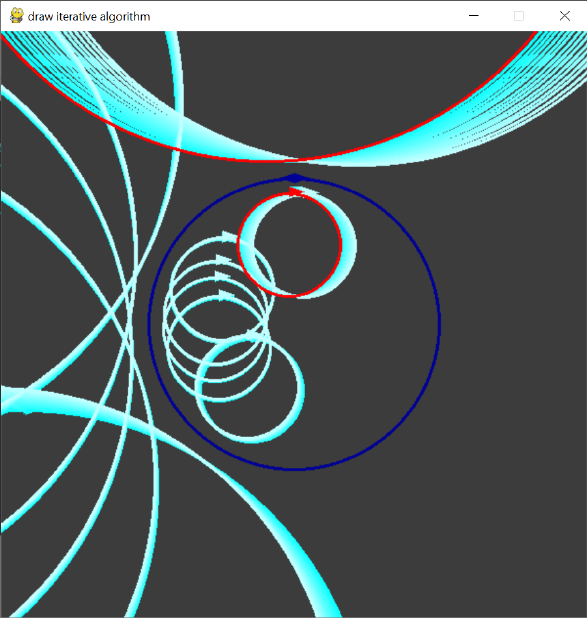}}
  \end{tabular}
  \caption{}
\end{table}

\

Consider the red oriented cycles in each case:
\begin{itemizedot}
  \item In case 1, we see that the \tmcolor{blue}{dark blue} circle meets the
  \tmcolor{red}{red cycle} at two distinct points. The two points where the
  \tmcolor{blue}{dark blue} and \tmcolor{red}{red} circles meet correspond to
  the \tmtextit{real eigenvectors}. Since in the two-circle case, we have
  $\det (M) \geq 0$, the eigenvalues of $M$ have the same sign. Therefore, $M$
  is \tmtextit{diagonalisable with real eigenvalues of the same sign}.
  
  \item In case 2, we see that the the \tmcolor{blue}{dark blue} circle meets
  the \tmcolor{red}{red cycle} at only one point. This one point is the sole
  real eigenvector of $M$. Therefore \tmtextit{$M$ is non-diagonalisable, and
  its sole eigenvalue is real}.
  
  \item In case 3, we see that the red cycle does not meet the dark blue
  circle. This means that $M$ has no real eigenvectors. But it must have
  \tmtextit{some} eigenvector. This eigenvector is complex, and has a complex
  eigenvalue $\lambda$. Since $\overline{\lambda}$ is also an eigenvalue, we
  have that \tmtextit{$M$ is diagonalisable with complex eigenvalues.}
\end{itemizedot}
These are 3 of the 4 possible Jordan Normal Forms of a $2 \times 2$ real
matrix. The remaining case has the eigenvalues of $M$ be real and with
opposite signs, which we will investigate later.

Observe that each \tmcolor{cyan}{light-blue} circle is obtained from the
\tmcolor{red}{red} circle by a rotation around the point in the plane $(0,
0)$. This means that each iteration is orthogonally similar to the preceeding
iterations. This is easily explained by the definition of the QR algorithm.

Observe the white circle(s): Clearly, in the first two cases, the algorithm
converges, and in the third, it doesn't. But notice that -- assuming
convergence -- the white circle passes through the point $(- 1, 0)$. A circle
passing through $(- 1, 0)$ corresponds to an upper triangular matrix. The QR
algorithm, when it converges, converges to an upper triangular matrix. The
\tmtextit{very objective} of the QR algorithm is to successively reduce the
angle between a cycle and the point $(- 1, 0)$.

Case 2 (the non-diagonalisable case) converges more slowly than case 1, but
still converges. There is a transition from cases 1, to 2, to 3, where
convergence gets slower until it stops happening.

We now explore behaviour near limits:

\begin{table}[h]
  \begin{tabular}{ll}
    4 & 5\\
    \raisebox{0.0\height}{\includegraphics[width=5.9119277187459cm,height=6.21802439984258cm]{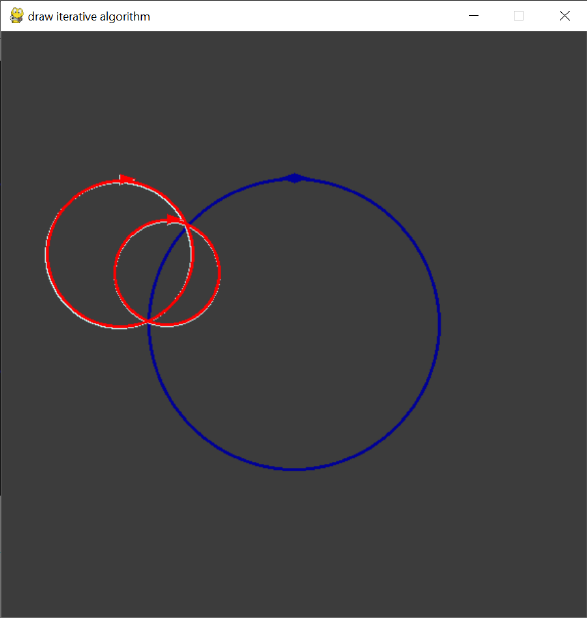}}
    &
    \raisebox{0.0\height}{\includegraphics[width=5.9119277187459cm,height=6.21802439984258cm]{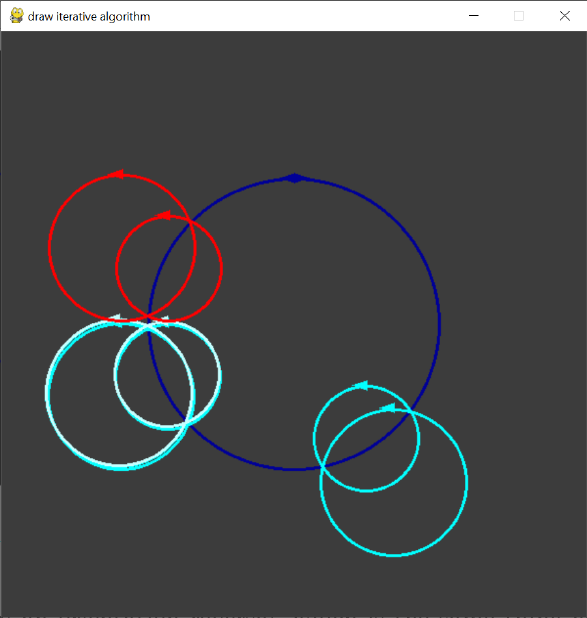}}
  \end{tabular}
  \caption{}
\end{table}

\

Case 5 is obtained from case 4 by an orientation reversal, which is the same
thing as taking the matrix inverse. Case 4 is an attractive fixed point; small
perturbations of it get mapped closer to the fixed point. Case 5 is a
repulsive fixed point; small perturbations are mapped away from the fixed
point. \tmtextit{We see that the matrix inverse of an attractive fixed point
is a repulsive fixed point.}

Without explaining too much about this point (see
{\cite{leite2013dynamics,goodshift2022}}: The existence of non-attractive
fixed points is a consequence of the fact that an iteration of $(Q R = M)
\mapsto R Q$ is continuous in $M$. It is a topological fact as opposed to a
quantitative fact. In order for a variant of QR to be practical, it must
iterate a \tmtextit{discontinuous} map. Observe that Wilkinson shifts, for
instance, introduce the necessary discontinuity.

Near the identity matrix, we see case 6:

\begin{table}[h]
  \begin{tabular}{ll}
    6 & 7\\
    \raisebox{0.0\height}{\includegraphics[width=5.9119277187459cm,height=6.21802439984258cm]{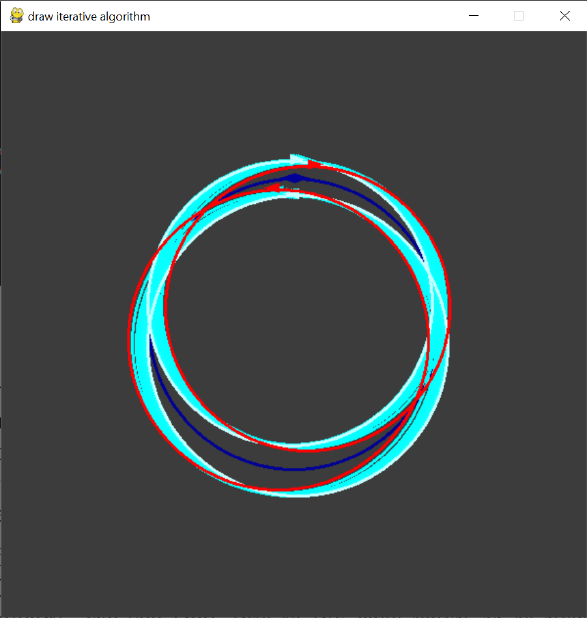}}
    &
    \raisebox{0.0\height}{\includegraphics[width=5.9119277187459cm,height=6.21802439984258cm]{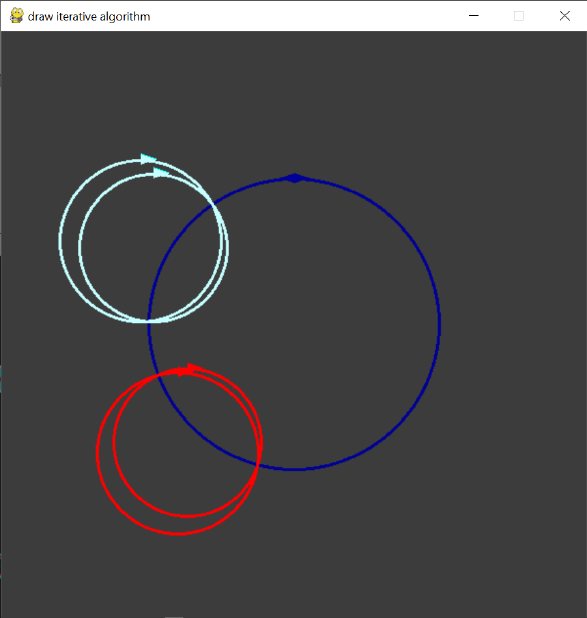}}
  \end{tabular}
  \caption{}
\end{table}

Case 6 has the two \tmcolor{red}{red cycles} be very close to the
\tmcolor{blue}{dark blue circle} (barely visible). This means they are close
to the identity matrix. In this case, the algorithm converges very slowly,
which can be seen from the large amounts of \tmcolor{cyan}{light blue}. This
behaviour happens because near an eigenvalue clash, the eigenspaces are
unstable. For such a matrix, while we can estimate its eigenvalues accurately
using the Gershgorin circle theorem, it is impossible to accurately determine
its eigenvectors.

The complete opposite of case 6 is case 7. In this case, the two red cycles
$\{ C, \tmop{ORI} (C) \}$ are almost the same as each other. They are also
nearly orthogonal to the dark blue circle, making them approximately
hyperbolic spears. We see that \tmtextit{this is the case when $M$ is nearly
singular}. Convergence happens very quickly here.

Case 6 happens when $\frac{\lambda_1}{\lambda_2} \approx 1$, and its behaviour
is thanks to fundamental computability considerations (finding eigenvectors
for a matrix whose eigenvectors have unknown algebraic multiplicities is
impossible), and case 7 happens when $\frac{\lambda_1}{\lambda_2} \approx 0$,
and its behaviour is the opposite to case 6.

The objective of \tmtextit{shifting} an $n \times n$ matrix is to
heuristically find a constant $\lambda$ such that $M - \lambda I$ is
approximately singular. This speeds up convergence. This is motivated by
moving away from case 6 to case 7.

In the complex eigenvalue case, the algorithm displays near periodicity,
except when $M$ is near to an orthogonal matrix:

\begin{table}[h]
  \begin{tabular}{ll}
    8 & 9\\
    \raisebox{0.0\height}{\includegraphics[width=5.9119277187459cm,height=6.21802439984258cm]{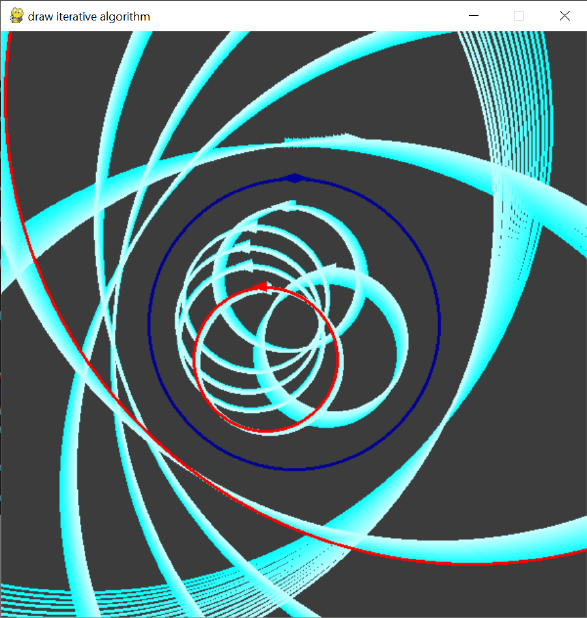}}
    &
    \raisebox{0.0\height}{\includegraphics[width=5.9119277187459cm,height=6.21802439984258cm]{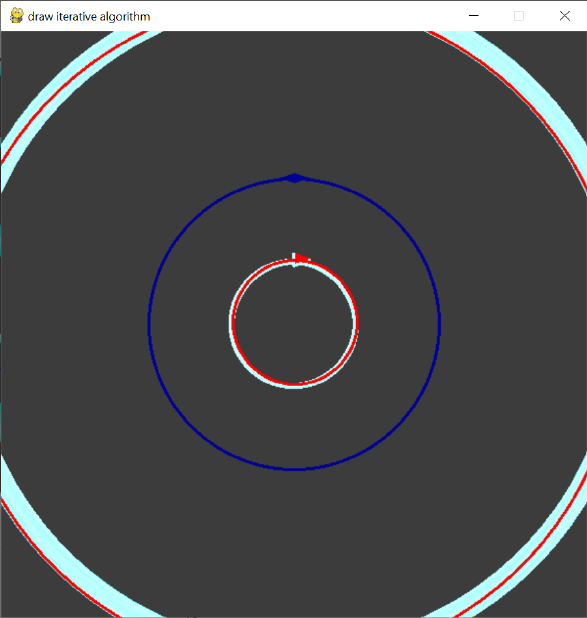}}
  \end{tabular}
  \caption{}
\end{table}

If the eigenvalues are $r e^{\pm i 2 \pi m / n}$ for $m / n$ in simplest form,
then the periodicity of the iterates is $n$. If the eigenvalues are instead $r
e^{\pm i 2 \pi x}$ for $x$ irrational, then the approximate periodicity is
given by the rational approximations of $x$.

Case 9 involves a near-orthogonal matrix. This is nearly a fixed point.

\subsection{The naive QR algorithm in the negative determinant case}

The user can rotate one of the arrow heads, and in so doing, the other arrow
head. The user can specify the pair of points where the hyperbolic line
representing the matrix meets the unit circle.

\tmtextit{Note: In the negative determinant case, we post-multiply Q and
pre-multiply $R$ by $\left(\begin{tabular}{cc}
  -1 & 0\\
  0 & 1
\end{tabular}\right)$ in order to make the dynamics a bit simpler. A user of
this visualisation may decide not to do this, and investigate.}

Below, we see some negative determinant cases. In these cases, the algorithm
converges to one limit point when $\tmop{trace} (M) < 0$ (case 10), and a
different one when $\tmop{trace} (M) > 0$ (case 12). When $\tmop{trace} (M) =
0$ (case 11), the algorithm oscillates between two points, both different from
the limit points in the other cases; one of the limit points is the starting
one; therefore, for small values of $\tmop{trace} (M)$, convergence is slow.

\begin{table}[h]
  \begin{tabular}{lll}
    10 & 11 & 12\\
    \raisebox{0.0\height}{\includegraphics[width=5.9119277187459cm,height=6.21802439984258cm]{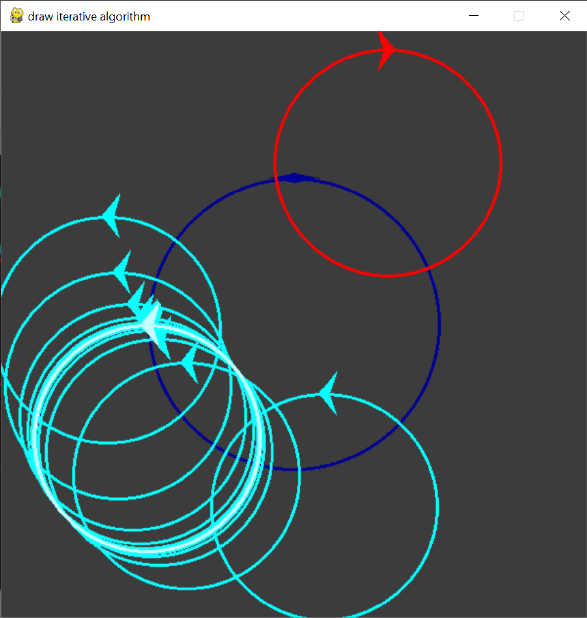}}
    &
    \raisebox{0.0\height}{\includegraphics[width=5.9119277187459cm,height=6.21802439984258cm]{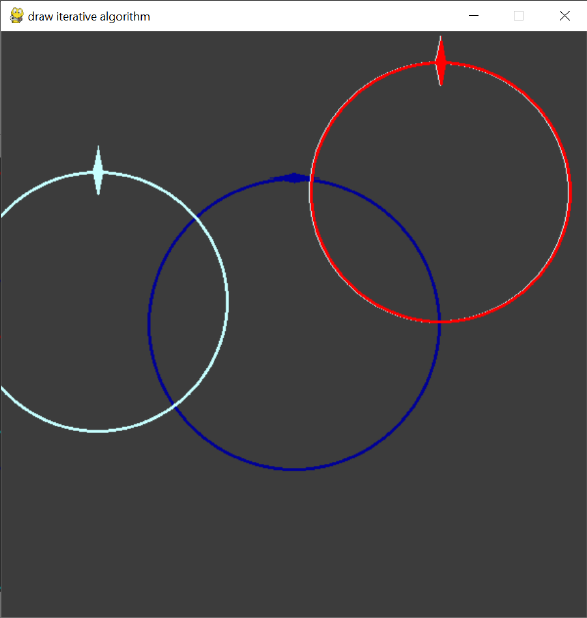}}
    &
    \raisebox{0.0\height}{\includegraphics[width=5.9119277187459cm,height=6.21802439984258cm]{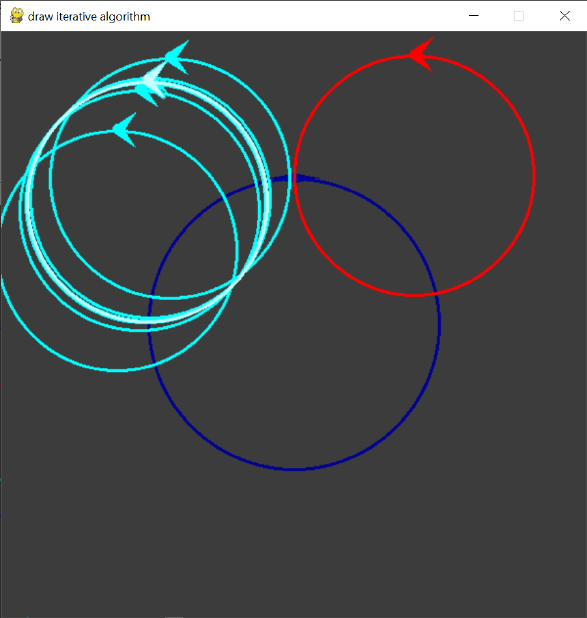}}
  \end{tabular}
  \caption{}
\end{table}

\section{How v2 works (Lie Sphere Geometry)}

\subsection{Lie Sphere Geometry vs Euclidean vs Moebius}

Here, we will introduce Lie Sphere Geometry {\cite{Benz2012}}, and use it to
verify the above claims. It should be possible to reproduce the visualisations
using this.

Lie Sphere Geometry (which we will abbreviate to LSG) is a generalisation of
Moebius geometry, which is in turn a generalisation of Euclidean geometry.

\subsection{Euclidean geometry}

Euclidean geometry is the study of objects called \tmtextit{cycles}, which are
either points, non-point circles, or lines. There is a group which acts on
these called the Euclidean group, which is denoted $E (n)$ in dimension $n$.
It is generated by the Euclidean reflections.

\subsection{Moebius geometry}

We turn to Moebius geometry. Synthetically (that is, without reference to
coordinates) Moebius geometry is the study of the same cycles as in Euclidean
geometry, except with an additional point at infinity. Euclidean motions
generalise so as to fix the point at infinity. The real difference lies in the
group, which consists of all conformal transformations, i.e. those which leave
angles between cycles the same.

To perform Moebius Geometry analytically, it is helpful to provide coordinate
systems. The points can be given homogeneous coordinates $[z : w]$, where $z$
and $w$ are both complex numbers. Thus, the points are those of the projective
line over the complex numbers, $\mathbb{\tmop{CP}}^1$. Clearly, the point with
homogeneous coordinate $[z : w]$ can be represented by the vector $(z, w)^T$.

The Moebius transformations form a group isomorphic to $\tmop{PGL}_2
(\mathbb{C})$. This group consists of the invertible $2 \times 2$ complex
matrices, modulo non-zero scalar factor.

The cycles in Moebius geometry can be coordinatised using \tmtextit{indefinite
$2 \times 2$ Hermitian matrices}, modulo real scalar factor. A point $p$ lies
on cycle $C$ iff $p^{\ast} C p = 0$, where $p^{\ast}$ denotes the
conjugate-transpose of $p$.

\subsection{Lie Sphere Geometry}

We now describe LSG. In LSG, the space consists of \tmtextit{oriented cycles}.
These are either unoriented points, oriented non-point circles, or oriented
lines. The transformation group consists of all those bijective maps which
preserve oriented tangency. The transformation contains all the Moebius
transformations and more; the non-Moebius transformations are not relevant to
our discussion.

The space of oriented cycles can be coordinatised using some of the points of
$\mathbb{\tmop{RP}}^4$. Not every point in $\mathbb{\tmop{RP}}^4$ represents
an oriented cycle. Rather, we define the \tmtextit{Lie quadric} $\mathbb{L}$,
whose points are in 1-to-1 correspondence with oriented cycles in the plane,
by $- x_1^2 - x_2^2 + x_3^2 + x_4 x_5 = 0$. Additionally, we define a method
for checking whether the oriented cycles represented by two points $p$ and $q$
of $\mathbb{L}$ are orientedly tangent. The method is to use the bilinear form
$B$ corresponding to the quadratic form which defines the Lie quadric
$\mathbb{L}$, and check whether $B (p, q) = 0$. We prefer the notation $p
\cdot q$ to $B (p, q)$. The transformation group is thus isomorphic to the
projective linear group $\tmop{PO}_{2, 3} (\mathbb{R})$.

The precise correspondence between points on $\mathbb{L}$ and oriented cycles
is as follows: A point $p \in \mathbb{L}$ either has:
\begin{enumerateroman}
  \item $p = [x_1 : x_2 : x_3 : x_4 : 1]$,
  
  \item or $p = [x_1 : x_2 : 1 : x_4 : 0]$,
  
  \item or $p = [0 : 0 : 0 : 1 : 0]$,
\end{enumerateroman}
and these correspond respectively to:
\begin{enumerateroman}
  \item An oriented circle with centre point $(x_1, x_2)$, radius $| x_3 |$,
  and orientation \tmtextit{inside the disk} if $\tmop{sign} (x_3) = - 1$, or
  \tmtextit{outside the disk} if $\tmop{sign} (x_3) = 1$.
  
  \item The oriented line given by the half-plane $x_1 x + x_2 y - 2 x_4 \geq
  0$.
  
  \item The point at infinity.
\end{enumerateroman}
Using this correspondence, one can then check that indeed $p \cdot q = 0$ iff
$p$ and $q$ are orientedly tangent.

\subsection{Applying Moebius transformations to oriented cycles}

An oriented cycle $C = [- b : - c : \ast : d : a]$ corresponds to the
unoriented cycle whose Hermitian matrix representation is $H =
\left(\begin{array}{cc}
  a & - b + i c\\
  - b - i c & d
\end{array}\right)$. The term hidden by $\ast$ is one of $\pm \sqrt{- \det
(H)}$. This correspondence enables one to evaluate Moebius transformations on
oriented cycles.

\subsection{Matrix as a bundle of hyperbolic spears}

It is easier initially to represent the projective real line as the $x$-axis
rather than as the unit circle. In this case, a nonnegative-determinant $2
\times 2$ real matrix is depicted as a pair $\{ C, \tmop{RI} (C) \}$, where
$C$ is any oriented cycle, and $\tmop{RI} (C)$ denotes the result of
reflecting $C$ through the $x$-axis and then reversing its orientation. This
looks as follows:

\raisebox{0.0\height}{\includegraphics[width=5.9119277187459cm,height=6.21802439984258cm]{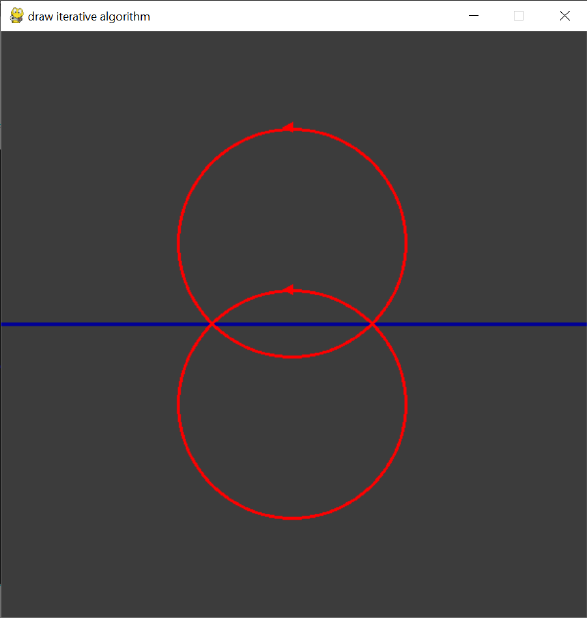}}

This is justified as follows: Every $2 \times 2$ real matrix $M$ acts on the
projective real line in a canonical way. We can represent the projective real
line as the $x$-axis within the ambient space $\mathbb{\tmop{RP}}^2$. We
define a \tmtextit{hyperbolic spear} (in this context) to be an oriented cycle
orthogonal to the $x$-axis. From any matrix $M$, we may produce the bundle
$\Gamma (M)$ of hyperbolic spears connecting each point $p$ on the $x$-axis to
$M p$. Depicting this for five values of $p$ may look like:

\raisebox{0.0\height}{\includegraphics[width=5.9119277187459cm,height=6.21802439984258cm]{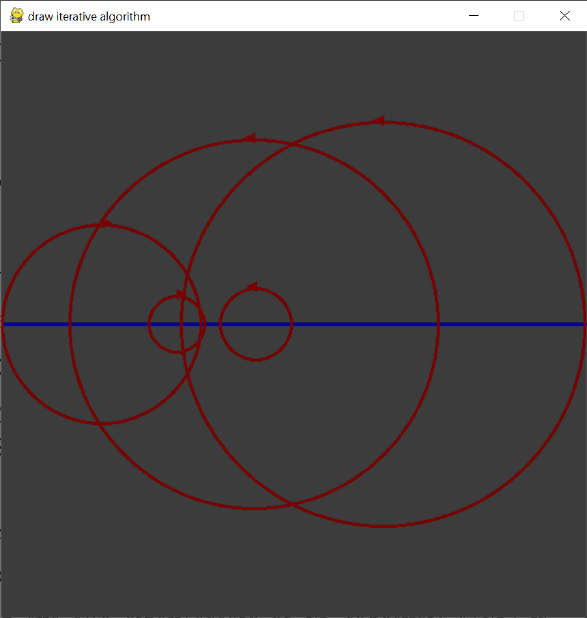}}

The dark red circles are the elements of $\Gamma (M)$.

The hyperbolic spear connecting $p = [x : y]$ to $M p$ where $M =
\left(\begin{array}{cc}
  a & b\\
  c & d
\end{array}\right)$ is given by $[- x (c x + d y) - y (a x + b y) : 0 : x (c x
+ d y) - y (a x + b y) : - 2 x (a x + b y) : - 2 y (c x + d y)]$.

We may reduce $\Gamma (M)$ to a simple figure in each case:

\tmtextit{In the positive determinant case:} Given matrix $M =
\left(\begin{array}{cc}
  a & b\\
  c & d
\end{array}\right)$, all the elements in $\Gamma (M)$ can be shown to be
orientedly tangent to $\{ C, \tmop{RI} (C) \}$ where $C = \left[ a - d : 2
\sqrt{a d - b c} : a + d : - 2 b : 2 c \right]$ and $\tmop{RI} (C) = \left[ a
- d : - 2 \sqrt{a d - b c} : a + d : - 2 b : 2 c \right]$. We see that this
naturally extends to the zero determinant case. We see also that this formula
\tmtextit{does not} work in the negative determinant case.

\tmtextit{In the negative determinant case:} Given a matrix $M =
\left(\begin{array}{cc}
  a & b\\
  c & d
\end{array}\right)$ with $\det (M) < 0$, all the elements of $\Gamma (M)$ can
be shown to meet some hyperbolic line $L (M)$ at angles $\pm \frac{\pi}{2} (1
+ \theta (M))$. The $L (M)$ is the unique hyperbolic line connecting the two
eigenvectors of $M$ together. $\theta (M)$ is $- 1 + \frac{2}{\pi} \arcsin
((\lambda_1 + \lambda_2) / (\lambda_1 - \lambda_2))$. These claims can be
verified by proving them in the case where $M$ is diagonal -- and $L (M)$ is
thus the $y$-axis -- and then generalising to non-diagonal $M = P D P^{- 1}$
(with $\det (D) < 0$) by treating $P$ as a Moebius transformation.

Finally, we apply a Moebius transform $z \mapsto \frac{1 - i z}{1 + i z}$ to
send the $x$-axis to the unit circle. This produces the visualisation we've
been using elsewhere in this paper.

\section{Prospects for a generalisation to $3 \times 3$ matrices}

\subsection{Brief introduction}

A possible criticism is that the above discussion is limited to $2 \times 2$
real matrices, where the QR algorithm is not strictly speaking used. In spite
of that, it is useful to convey intuitively how the naive, uncomplicated
algorithm behaves on the simplest inputs. The behaviour is complicated enough
that a visualisation ought to have value. Here, we suggest an extension to the
$3 \times 3$ case.

Some aspects of v2 can formally be generalised to $3 \times 3$ real matrices.
Once again, one tries to take an \tmtextit{envelope} of some oriented figures.
The problem is that the figures are now slightly more complicated, and the
envelope -- like already encountered in the $2 \times 2$ case -- might not
exist. In the $2 \times 2$ case, we know how to deal with this problem: If a
matrix has negative determinant, then we need to use a certain
\tmtextit{generalisation} of an envelope that we described above. What the
appropriate generalisation might be in the $3 \times 3$ case is not clear.

\subsection{Details}

Let $M$ be a $3 \times 3$ real matrix. $M$ acts in a natural way on the
projective plane $\mathbb{\tmop{RP}}^2$. While an element of
$\mathbb{\tmop{RP}}^2$ may be a point, it is actually more convenient to
choose it to be a\tmtextit{ projective line}, thanks to projective duality.
Embed $\mathbb{\tmop{RP}}^2$ as a plane $\Pi$ within the ambient projective
space $\mathbb{\tmop{RP}}^3$. Given a line $l$ on $\Pi$, $M$ maps it to $M
(l)$. Connect $l$ to $M (l)$ using an \tmtextit{oriented right circular cone}.
This oriented figure will necessarily exist and be unique, even when $M (l)$
is the unique line at infinity on $\Pi$. Finally, one must (if possible) take
the envelope of all of these oriented right circular cones, each one
corresponding to a different line $l$. The resulting envelope may or may not
exist, but some sort of generalisation of an envelope may be applicable (see
the discussion above).

\paragraph{Acknowledgements}I would like to thank Jan Techter, Gregory Gutin
and Abbas Edelat for some helpful discussions.


\begin{thebibliography}{1}
  \bibitem[1]{goodshift2022}Jess Banks, Jorge Garza-Vargas, and  Nikhil
  Srivastava. {\newblock}Global convergence of hessenberg shifted qr i:
  dynamics. {\newblock}2021.{\newblock}
  
  \bibitem[2]{Benz2012}Walter Benz. {\newblock}\tmtextit{Sphere Geometries of
  M{\"o}bius and Lie},  pages  93--174. {\newblock}Springer Basel, Basel,
  2012.{\newblock}
  
  \bibitem[3]{buurema1958geometric}Hendrik~Jan Buurema. {\newblock}\tmtextit{A
  geometric proof of convergence for the $QR$ method}.
  {\newblock}Rijksuniversiteit te Groningen, Groningen, 1970.
  {\newblock}Doctoral dissertation, University of Groningen.{\newblock}
  
  \bibitem[4]{leite2013dynamics}Ricardo~S.~Leite, Nicolau~C.~Saldanha, and 
  Carlos Tomei. {\newblock}Dynamics of the symmetric eigenvalue problem with
  shift strategies. {\newblock}\tmtextit{Int. Math. Res. Not. IMRN},
  (19):4382--4412, 2013.{\newblock}
  
  \bibitem[5]{watkins2008qr}David~S.~Watkins. {\newblock}The $QR$ algorithm
  revisited. {\newblock}\tmtextit{SIAM Rev.}, 50(1):133--145, 2008.{\newblock}
\end{thebibliography}
\end{document}